\documentclass[12pt]{amsart}

\begin{document}
\title[Moore-Penrose Inverse in $C^*$-algebras]{ Moore-Penrose inverse and doubly commuting elements in $C^*$-algebras}
\author[Enrico Boasso]{Enrico Boasso}

\begin{abstract} In this work it is proved that the Moore-Penrose inverse of the product
of $n$-doubly commuting regular $C^*$-algebra elements obeys the so-called reverse order law.
Conversely, conditions regarding the reverse order law of the Moore-Penrose
inverse are given in order to characterize when $n$-regular elements doubly commute. 
Furthermore, applications of the main results of this article to normal $C^*$-algebra
elements, to Hilbert space operators and to Calkin algebras will be considered. 
\end{abstract}
\dedicatory{This work is dedicated to the memory of Professor Angel R. Larotonda,\\
Pucho for all who knew him.}
\subjclass[2000]{46L05, 47A05}
\keywords{Generalized inverse, Moore-Penrose inverse, and doubly commuting 
elements in a $C^*$-algebra.}
\maketitle
\noindent \bf{1. Introduction}\rm \vskip.3cm

\indent Consider an unitary ring $A$. An element $a\in A$ will be said \it regular
\rm if it has a \it generalized inverse \rm in A, that is if there exists $b\in A$ such that
$$
a=aba.
$$  
A generalized inverse is also termed a \it pseudo inverse\rm .\par

\indent Note that if $a$ is a regular element of $A$ and $b$ is a generalized
inverse of $a$, then $p=ba$ and $q=ab$ are \it idempotents \rm of A, that is
$p=p^2$ and $q=q^2$.\par
 
\indent Given $a\in A$ a regular element, a generalized inverse $b$ of $a$
will be called \it normalized\rm, if $b$ is regular and $a$ is a pseudo 
inverse of $b$, equivalently if
$$
a=aba,\hskip1cm b=bab.
$$
Recall that if $b$ is a generalized inverse of $a$, then 
$c=bab$
is a normalized pseudo inverse of $a$.\par
\indent Next suppose that $a$ is a regular element and $b$ is a normalized
generalized inverse of $a$. In the presence of an involution $*\colon A\to A$, it is
possible to enquire if the idempotents $p$ and $q$ are \it self-adjoint\rm ,
that is whether $(ba)^*=ba$ and $(ab)^*=ab$. In this case $b$ is called the \it Moore-Penrose inverse \rm of $a$, see [16] where
this concept was introduced.\par
\indent In [10] it was proved that each regular element $a$ in a $C^*$-algebra
$A$ has a uniquely determined Moore-Penrose inverse. The
Moore-Penrose inverse of $a\in A$ will be denote by $a^{\dag}$. Therefore,
the Moore-Penrose inverse of a regular element $a\in A$ is the unique
element that satisfy the following equations:
$$
a=aa^{\dag}a,\hskip.5cm a^{\dag}=a^{\dag}aa^{\dag},\hskip.5cm (a^{\dag}a)^*=a^{\dag}a,\hskip.5cm
(aa^{\dag})^*=aa^{\dag}.
$$ 
\indent According to the uniqueness of the notion under consideration, if $a$ has a Moore-Penrose inverse,
then $a^*$ also has a Moore-Penrose inverse and
$$
(a^*)^{\dag}=(a^{\dag})^*.
$$
\indent Moreover, according to the above equations, if $a$ is a regular element,
then $a^{\dag}$ also is and 
$$
(a^{\dag})^{\dag}=a.
$$
\indent For other properties regarding Moore-Penrose inverses in 
$C^*$-algebras, see the works [10], [11], [14] and [16].\par 

\indent On the other hand, an $n$-tuple $a=(a_1,\dots ,a_n)$ of elements in a $C^*$-algebra A will
be called \it doubly commuting\rm, if $a_ia_j=a_ja_i$ and $a_ia_j^*=a_j^*a_i$,
for all $i$, $j=1,\ldots ,n$, $i\neq j$. For instance, necessary and sufficient for $(a,a)$
to be doubly commuting is that $a$ is a normal element in $A$. \par

\indent Doubly commuting operators have been  
studied in very different contexts, to mention only some of the most relevant works, 
see for example [1]-[4], [6]-[8] and [12].
In this article, doubly commuting tuples of regu-
lar $C^*$-algebra elements will be
consider. The main objective of this work consists in the study of the relationship
between such tuples and the Moore-Penrose inverse of the product of the 
elements in the tuple.\par

\indent In fact, given an $n$-tuple $a=(a_1,\dots ,a_n)$ of regular elements in a 
$C^*$-algebra A, if the $n$-tuple $a$ is doubly commuting, then $\prod_{i=1}^n a_i$
is regular and 
$$
(\prod_{i=1}^n a_i)^{\dag}\hbox{ }=\hbox{ }\prod_{i=0}^{n-1} a_{n-i}^{\dag}=
\hbox{ }\prod_{i=1}^n a_i^{\dag},
$$
in particular, $\prod_{i=1}^n a_i$ complies with the so-called reverse order law for
the Moore-Penrose inverse. Moreover, if $a$ is such a tuple, then $(a_i,a_j)$
and $(a_i^*,a_j)$ are doubly commuting pairs, where $i$, $j\in[\![ 1,n]\!]$,
$i\neq j$. Consequently, $a_ia_j$ and $a_i^*a_j$ are regular and
$$
(a_ia_j)^{\dag}=(a_ja_i)^{\dag}, \hskip1cm (a_i^*a_j)^{\dag}=(a_ja_i^*)^{\dag}.
$$ 
\indent Conversely, if $a=(a_1,\dots ,a_n)$ is an $n$-tuple of regular elements in a 
$C^*$-algebra A such that the above identities are verified by 
$(a_i,a_j)$ and $(a_i^*,a_j)$, 
for all $i$, $j\in[\![ 1,n]\!]$, $i\neq j$, then $a$ is a doubly commuting $n$-tuple
of elements in the $C^*$-algebra $A$. \par

\indent It is worth noticing that this characterization consists in an extension to the objects
under consideration of the sufficient conditon given by J.J. Koliha in [13; 2.13] for the product of two regular
elements to be Moore-Penrose invertible.\par
 
\indent In section 2 it will be proved the aforementioned characterization. 
In section three,
on the other hand, some applications will be developed. In fact, three sorts of
objects will be considered, namely,
tuples of commuting regular normal $C^*$-algebra 
elements, tuples of doubly commuting regular Hilbert space operators, and tuples 
of almost doubly commuting regular Hilbert space operators 
(that is doubly commuting tuples of regular elements in Calkin algebras).\par

\indent This article is dedicated to the memory of Professor Angel R. Larotonda,
who unfortunately and unexpectedly died on January 2th 2005. Although it is 
not necessary to comment Professor Larotonda's work 
as mathematician, for his scientific publications consist in a set of achievements
which speak for themselves, a few words
about the man deserve to be said. The author knew Professor Larotonda for 
more than twenty years. During uncountable conversations shared
with Professor Larotonda, this researcher always showed his condition of sensible, civilized 
and cultivated human being, three characteristics that seem to be far from being
widespread in this time and in any time.\par      
\vskip.3cm
\noindent \bf{Acknowledgements.} \rm The author wishes to express his indebtedness to the 
organizers of this volume, especially to Professor G. Corach, for have invited the
author to contribute with this homage.\par 

\vskip.3cm
\noindent \bf{2. Main Results}\rm \vskip.3cm

\indent In this section the relationship between the Moore-Penrose inverse and 
doubly commuting tuples of regular $C^*$-algebra elements will be studied. 
In fact, the characterization described in the previous section will be proved.
In first place,
a property of doubly commuting tuples is discussed.\par

\newtheorem*{rem2.1}{Remark 2.1}
\begin{rem2.1}\rm Let $a=(a_1,\ldots ,a_n)$ be an $n$-tuple of elements
of a $C^*$-algebra $A$. Consider $\pi\colon [\![ 1,n]\!]\to [\![ 1,n]\!]$ a permutation, and 
define an $n$-tuple $b=(b_1,\ldots ,b_n)$ in such a way that $b_j$ is  
either $a_{\pi(j)}$ or $a_{\pi(j)}^*$. Then, it is easy to prove that $a$ is doubly
commuting if and only if $b$ is.\par
\indent Furthermore, note that the following facts are equivalents:\par
\hskip.3cm i) $a=(a_1,\ldots ,a_n)$ is a doubly commuting 
$n$-tuple of elements of $A$,\par
\hskip.3cm ii) for each $i$, $j=1,\ldots ,n$, $i\neq j$,
$a_{i,j}=(a_i,a_j)$ is a pair of doubly commuting elements of $A$.\par
\end{rem2.1}  
\newtheorem*{prop2.2}{Proposition 2.2}
\begin{prop2.2} Let $a=(a_1,\ldots ,a_n)$ be an $n$-tuple of regular elem-
ents in 
a $C^*$-algebra $A$. Consider the $n$-tuple $a^{\dag}=(a_1^{\dag},\ldots ,a_n^{\dag})$.
Then, $a$ is doubly commuting if and only if $a^{\dag}$ is. \par
\end{prop2.2}
\begin{proof} According to Remark 2.1, it is enough to prove that
a pair of regular elements $(b,c)$ is doubly commuting if and
only if $(b^{\dag},c^{\dag})$ is.\par
\indent Suppose that $(b,c)$ is a doubly commuting pair of regular
elements of $A$.
Then, according to Theorem 5 of [10], 
$$
b^{\dag}c=cb^{\dag}.
$$
Moreover, since according to Remark 2.1 $(b,c^*)$ is a doubly commuting
pair, it is clear that 
$$
b^{\dag}c^*=c^*b^{\dag}.
$$
Consequently, $(b^{\dag},c)$ is a doubly commuting pair.\par
\indent However, according to Remark 2.1, $(c,b^{\dag})$ is a doubly 
commuting pair, and thanks to what has been proved,
$(c^{\dag},b^{\dag})$ is a doubly com-
muting pair. Therefore, according again
to Remark 2.1, $(b^{\dag},c^{\dag})$
is a doubly commuting pair.\par
\indent Conversely, if $(b^{\dag},c^{\dag})$ is a doubly commuting pair,
since $(b^{\dag})^{\dag}=b$ and $(c^{\dag})^{\dag}=c$, according to the first part of the proof,
$(b,c)$ is a doubly commuting pair 
\end{proof}

\indent Note that in [13] Theorem 5 of [10] was proved using the Drazin inverse,
see [13; 2.12].
Therefore, Proposition 2.2 can also be derived using the Drazin inverse,
see [13; 2.12] and [13; 2.13].\par

\newtheorem*{rem2.3}{Remark 2.3}
\begin{rem2.3}\rm Let $a=(a_1,\ldots ,a_n)$ be an $n$-tuple of regular elements in 
a $C^*$-algebra $A$. Consider $\pi\colon [\![ 1,n]\!]\to [\![1,n]\!]$ a permutation, and define an
$n$-tuple $b=(b_1,\ldots ,b_n)$ in the following way. Given $j=1,\ldots ,n$,
$b_j$ is either $a_{\pi(j)}$, $a_{\pi(j)}^*$, $a_{\pi(j)}^{\dag}$ or $(a_{\pi(j)}^{\dag})^*$. 
Then, according to 
Remark 2.1  and to Proposition 2.2, $a$ is doubly commuting if and only if
$b$ is. \par
\indent Furthermore, according to Remark 2.1 and to Proposition 2.2, the following facts are equivalent:\par
\hskip.3cm i) $a$ is an $n$-tuple of doubly commuting regular elements of $A$,\par
\hskip.3cm ii) for each $i$, $j=1,\ldots ,n$, $i\neq j$, $(b_i,b_j)$ is a pair of doubly 
commuting regular elements of $A$. \par 
\end{rem2.3}
\indent Next follows the first part of our characterization. In fact, in the following theorem
it will be proved that the product of the elements in a doubly commuting tuple of regular elements satisfy the 
so-called reverse order law for the Moore-Penrose inverse.\par
\newtheorem*{theo2.4}{Theorem 2.4}
\begin{theo2.4}  Let $a=(a_1,\ldots ,a_n)$ be an $n$-tuple of doubly commuting  
regular elements of 
a $C^*$-algebra $A$. Then, $\prod_{i=1}^n a_i$ is regular and
$$
(\prod_{i=1}^n a_i)^{\dag}=\prod_{i=0}^{n-1} a_{n-i}^{\dag}=\prod_{i=1}^{n} a_i^{\dag}.
$$
\end{theo2.4}
\begin{proof} Consider $b$ and $c$ two regular elements of $A$ such
that the pair $(b,c)$ is doubly commuting. According to [13; 2.13],
$bc$ is regular and
$$
(bc)^{\dag}=c^{\dag}b^{\dag}.
$$
\indent Next consider $a=(a_1,\ldots ,a_n,a_{n+1})$ a doubly commuting
tuple of regular elements of $A$. Suppose that $\prod_{i=1}^na_i$
is regular and 
$$
(\prod_{i=1}^na_i)^{\dag}=\prod_{i=0}^{n-1}a_{n-i}^{\dag}. 
$$
\indent In order to conclude the proof of the reverse order law for the
Moore- Penrose inverse, it is enough to prove that
$\prod_{i=1}^{n+1}a_i=(\prod_{i=1}^na_i)a_{n+1}$ is regular and
$$
(\prod_{i=1}^{n+1}a_i)^{\dag}=a_{n+1}^{\dag}(\prod_{i=0}^{n-1}a_{n-i})^{\dag}.
$$
\indent Consider $d=\prod_{i=1}^na_i$. According to Remark 2.1, $(a_{n+1},a_j)$ is a doubly 
commuting pair for $j=1,\ldots ,n$, which implies that $(d,a_{n+1})$ is a doubly 
commuting pair of regular elements. Therefore, according to the part of the theorem that has been
proved, $da_{n+1}=\prod_{i=1}^{n+1}a_i$
is regular, and
$$
(\prod_{i=1}^{n+1}a_i)^{\dag} =(da_{n+1})^{\dag}=a_{n+1}^{\dag}d^{\dag}
=\prod_{i=0}^na_{n+1-i}^{\dag}.
$$
\indent The last identity is a consequence of  [13; 2.13] or Proposition 2.2.\par
\end{proof}

\newtheorem*{rem2.5}{Remark 2.5}
\begin{rem2.5}\rm Let $a=(a_1,\ldots ,a_n)$ be an $n$-tuple of doubly commuting regular
elements in a $C^*$-algebra $A$. Consider $\pi\colon[\![1,n]\!]\to [\![1,n]\!]$ a permutation
and consider an $n$-tuple $b=(b_1\ldots ,b_n)$, where $b_j$ is either $a_{\pi(j)}$,
$a_{\pi(j)}^*$, $a_{\pi(j)}^{\dag}$ or $(a_{\pi(j)}^{\dag})^*$, $j=1,\ldots ,n$. Then, $b$ is an $n$-tuple of doubly commuting regular 
elements of $A$. Consequently, according to [13; 2.13] or Proposition 2.2 and to Theorem 2.4, $\prod_{i=1}^n b_i$
is regular and
$$
(\prod_{i=1}^n b_i)^{\dag}=\prod_{i=0}^{n-1} b_{n-i}^{\dag}=\prod_{i=1}^{n} b_i^{\dag}.
$$ 
\indent Next consider two permutations $\pi, \sigma\colon[\![1,n]\!]\to [\![1,n]\!]$.
Let $\tau\colon[\![1,n]\!]\to [\![1,n]\!]$ be the permutation $\tau=\pi^{-1}\sigma$,
that is $\pi\tau=\sigma$. Next associate to $\pi$ an $n$-tuple $b$ as in the
previous paragraph. Then, it is possible to associate to $\sigma$ an $n$-tuple $c$
such that $b_{\tau (j)}=c_j$, for $j=1,\ldots ,n$. Consequently, according
to [13; 2.13] or Proposition 2.2 and to Theorem 2.4,
$$
(\prod_{i=1}^n b_i)^{\dag}=\prod_{i=1}^n b_i^{\dag}=
\prod_{i=1}^n c_i^{\dag}=(\prod_{i=1}^n c_i)^{\dag}.
$$
\indent For instance, if $a=(b,c)$ is a pair of doubly commuting regular elements 
of the $C^*$-algebra $A$, the following identities hold.
\begin{align*}
(ab)^{\dag}&=b^{\dag}a^{\dag}=a^{\dag}b^{\dag}=(ba)^{\dag},\\
(a^*b)^{\dag}&=b^{\dag}(a^{\dag})^*=(a^{\dag})^*b^{\dag}=(ba^*)^{\dag},\\
(ab^*)^{\dag}&=(b^{\dag})^*a^{\dag}=a^{\dag}(b^{\dag})^*=(b^*a)^{\dag},\\
(a^*b^*)^{\dag}&=(b^{\dag})^*(a^{\dag})^*=(a^{\dag})^*(b^{\dag})^*=(b^*a^*)^{\dag}.
\end{align*} 
\indent Consequently,
\begin{align*}
(ab)^{\dag}&=(ba)^{\dag}=((a^*b^*)^{\dag})^*=((b^*a^*)^{\dag})^*,\\
(a^*b)^{\dag}&=(ba^*)^{\dag}=((ab^*)^{\dag})^*=((b^*a)^{\dag})^*.
\end{align*}
\indent Finally, recall that according to Remark 2.1, $a=(a_1,\ldots ,a_n)$
is a doubly commuting $n$-tuple of regular elements of $A$ if and only if
$(a_i,a_j)$ is a doubly commuting pair of regular elements of $A$, for
$i$, $j=1,\ldots ,n$, $i\neq j$. Thererfore, for each such pair the above identities
hold.\par
\end{rem2.5}
\indent Next follows the converse of Theorem 2.4.\par

\newtheorem*{theo2.6}{Theorem 2.6} 
\begin{theo2.6}Consider $a=(a_1,\ldots ,a_n)$ an $n$-tuple 
of regular elem-
ents in a $C^*$-algebra $A$. Suppose that for all $i$, $j=1,\ldots ,n$, 
$i\neq j$, $a_ia_j$, $a_ja_i$, $a_i^*a_j$ and $a_ja_i^*$ are regular and
comply with the following identities:
$$
(a_ia_j)^{\dag} =(a_ja_i)^{\dag},\hskip1cm (a_i^*a_j)^{\dag}=(a_ja_i^*)^{\dag}.
$$
Then $a$ is a doubly commuting $n$-tuple of regular elements of $A$.
\end{theo2.6}
\begin{proof}
\indent According to Remark 2.1, it is enough to prove that $(a_i,a_j)$
is doubly commuting for all $i$, $j=1,\ldots ,n$, $i\neq j$. However,
since $a_ia_j$ and  $a_ja_i$ are regular and 
$(a_ia_j)^{\dag} =(a_ja_i)^{\dag}$, it is clear that $a_ia_j =a_ja_i$. Similarly,
$a_i^*a_j=a_ja_i^*$. Therefore,  $(a_i,a_j)$ is a doubly commuting pair,
$i$, $j=1,\ldots ,n$, $i\neq j$.\par
\end{proof}

\newtheorem*{Rem2.7}{Remark 2.7}
\begin{Rem2.7} \rm In [13] sufficient conditions for a product 
to be Moore-Penrose inversible were given. Actually, thanks to the 
Theorem 2.6 it is now possible to state the following characterization.
Let $a$ and $b$ be two regular elements in a $C^*$-algebra A. If $(a,b)$
is a doubly commuting pair, then $ab$, $ba$, $a^*b$ and $ab^*$ are
regular and
$$
(ab)^{\dag} =(ba)^{\dag},\hskip1cm (a^*b)^{\dag}=(ba^*)^{\dag}.
$$
\indent On the other hand, if $a$ and $b$ are two regular elements 
such that $ab$, $ba$, $a^*b$ and $ab^*$ are
regular and they comply with the above identities, then $(a,b)$ 
is a doubly commuting pair. Furthermore, in this case $ab$ and $a^*b$ comply the
reverse order law for the Moore-Penrose inverse.
\par
\end{Rem2.7}

\noindent \bf{3. Some applications}\rm \vskip.3cm
\indent In this section several applications of the main results of this
work will be considered. In first place, commuting tuples 
of normal elements in a $C^*$-algebra will be studied.\par

\newtheorem*{Rem3.1}{Remark 3.1}
\begin{Rem3.1} \rm Let $a=(a_1,\ldots ,a_n)$ be a doubly commuting $n$-tuple
of elements in a $C^*$-algebra $A$. Let $\alpha =(m_1,\ldots ,m_n)$ be an $n$-tuple
of nonnegative entire numbers and define the $n$-tuple 
$$
a_{\alpha}=(a_1^{m_1},\ldots ,a_n^{m_n}).
$$
Then, it is not difficult to prove that $a_{\alpha}$ is a doubly commuting 
$n$-tuple of elements of $A$.\par

\indent On the other hand, if $b\in A$ is a regular normal element of $A$, then $b^n$ 
is regular
and $(b^n)^{\dag}=(b^{\dag})^n$, for all $n\in\Bbb{N}$. \par 

\indent In fact, since $b$ is regular, $b_n=(b,b,\ldots ,b)$
($n$-times) is a doubly commuting $n$-tuple of regular elements of $A$.
Therefore, according to Theorem 2.4, $b^n$ is regular and $(b^n)^{\dag}= (b^{\dag})^n$.\par 
\end{Rem3.1}
\indent In the following theorem the characterization of the previous section will
be applied to commuting tuples of regular normal elements. This result is an
extension  of [13; 2.14] to the case under consideration.\par 
\newtheorem*{theo3.2}{Theorem 3.2}
\begin{theo3.2}Let $a=(a_1,\ldots ,a_n)$ be a commuting $n$-tuple of 
regular normal elements in a $C^*$-algebra $A$. Let $\alpha =(m_1,\ldots ,m_n)$ be an $n$-tuple
of nonnegative entire numbers. Then $\prod_{i=1}^n a_i^{m_i}$ is regular
and 
$$
(\prod_{i=1}^n a_i^{m_i})^{\dag}=\prod_{i=1}^n (a_i^{\dag})^{m_i}.
$$
\end{theo3.2}
\begin{proof} First of all recall that according to the Flugede-Putman Theorem, 
[5; IX, 6.7] or [9; 9.6.7], $a=(a_1,\ldots ,a_n)$ is a doubly
commuting tuple.\par
Next consider, as in Remark 3.1, the $n$-tuple
$$
a_{\alpha}=(a_1^{m_1},\ldots ,a_n^{m_n}).
$$
\indent Then, according to Remark 3.1, $a_{\alpha}$ is a doubly commuting
tuple of regular elements of $A$. Consequently, according to Theorem 2.4, 
Remark 2.5 and Remark 3.1 again, the proof of the Theorem is concluded.\par  
\end{proof}

\newtheorem*{Rem3.3}{Remark 3.3}
\begin{Rem3.3}\rm Let $a=(a_1,\ldots ,a_n)$ be a commuting
$n$-tuple of regu-
lar normal elements in 
a $C^*$-algebra $A$. Consider, as in Remark 2.3, 
$\pi\colon [\![ 1,n]\!]\to [\![1,n]\!]$ a permutation, and define an
$n$-tuple $b=(b_1,\ldots ,b_n)$, where given $j=1,\ldots ,n$,
$b_j$ is either $a_{\pi(j)}$, $a_{\pi(j)}^*$, $a_{\pi(j)}^{\dag}$ or $(a_{\pi(j)}^{\dag})^*$.
Next consider the $n$-tuple $b_{\beta}=(b_1^{m_1},\ldots ,b_n^{m_n})$,
where  $\beta =(m_1,\ldots ,m_n)$ is an $n$-tuple of nonnegative integers.
Then, according to Proposition 2.2, Remark 3.1, Theorem 3.2 and Theorem 10 of [10],
$b_{\beta}$ is an $n$-tuple of doubly commuting regular normal elements of
$A$. Therefore, according to Theorem 3.2, $\prod_{i=1}^n b_i^{n_i}$ is regular
and 
$$
(\prod_{i=1}^n b_i^{n_i})^{\dag}=\prod_{i=1}^n (b_i^{\dag})^{n_i}.
$$
\indent For example, if $(c,d)$ is a pair of commuting regular 
normal elements in a $C^*$-algebra $A$, then 
\begin{align*}
& a^kb^l,&  &(a^k)^*b^l,&  &a^k(b^l)^*,&  &(a^k)^*(b^l)^*,\\
 &(a^k)^{\dag}b^l,&  &a^k(b^l)^{\dag},& &(a^k)^{\dag}(b^l)^{\dag},& &((a^k)^{\dag})^*((b^k)^{\dag})^*,\\
&((a^k)^{\dag})^*b^l,& &(a^k)^{\dag}(b^l)^*,& &((a^k)^{\dag})^*(b^l)^*,&
&a^k((b^l)^{\dag})^*,\\
&(a^k)^*(b^l)^{\dag},& &(a^k)^*((b^l)^{\dag})^*,&
&((a^k)^{\dag})^*(b^l)^{\dag},& &(a^k)^{\dag}((b^l)^{\dag})^*
\end{align*}
\noindent are regular elements of $A$, for $k$ and $l\in\Bbb{N}$. Furthermore,
their Moore-Pensrose inverses can be calculated according to
the first part of the present Remark.\par
\end{Rem3.3}

\indent Next the main results of the present work will be applied to $n$-tuples of
regular Hilbert space operators.\par

\indent Let $H$ be a Hilbert space and consider $A=L(H)$, the $C^*$-algebra of all
bounded and linear maps defined in $H$. Recall that $T\in L(H)$ is regular
as an operator if and only if $T$ is a regular element of $A$. Moreover,
necessary and sufficient for $T\in L(H)$ to be a regular operator
is the fact that the range of $T$, $R(T)$, is a closed subspace of $H$, see
for example [9; 3.8]. Note that in this case the Moore-Penrose inverse
can be described in a direct way.\par
\indent In fact, consider $T\in L(H)$ a bounded  Hilbert space operator with 
close range. Define $S\in L(H)$ as follows:
\begin{align*}
&S\colon R(T)\to N(T)^{\perp}, \hskip.8cm S\mid_{R(T)}\equiv \tilde{T}^{-1},\\
&S\colon R(T)^{\perp}\to H,\hskip1.4cm S\mid_{R(T)^{\perp}}\equiv 0,
\end{align*}
where $\tilde{T}=T\mid_{N(T)^{\perp}}^{R(T)}\colon N(T)^{\perp}\to R(T)$,
that is the restriction to $N(T)^{\perp}$ of $T$. Then, it is not diffuclt to
prove that
$$
T=TST,\hskip.3cm S=STS,\hskip.3cm (TS)^*=TS,\hskip.3cm (ST)^*=ST,
$$
that is, $S$ is the Moore-Penrose inverse of $T$.\par
 
\indent On the other hand, an $n$-tuple of continuous linear maps
defined in $H$, $T=(T_1,\ldots ,T_n)$, is doubly commuting as operators if and only 
if it is doubly commuting as elements of $A=L(H)$. In the following theorems the 
relationship between doubly commuting tuples of regular operators and
the Moore-Penrose inverse will be studied.\par

\newtheorem*{theo3.4}{Theorem 3.4}
\begin{theo3.4} Let $T=(T_1,\ldots ,T_n)$ be an $n$-tuple of 
regular Hilbert space operators. Then, if the $n$-tuple $T$ is doubly commuting,
$\prod_{i=1}^nT_i$ is a regular operator and
$$
(\prod_{i=1}^nT_i)^{\dag} =\prod_{i=0}^{n-1}T_{n-i}^{\dag}
=\prod_{i=1}^nT_i^{\dag}.
$$ 
\indent Conversely, if for all $i$, $j=1,\ldots ,n$, 
$i\neq j$, $T_iT_j$, $T_jT_i$, $T_i^*T_j$ and $T_jT_i^*$ are regular operators
which comply with the following identities:
$$
(T_iT_j)^{\dag} =(T_jT_i)^{\dag},\hskip1cm (T_i^*T_j)^{\dag}=(T_jT_i^*)^{\dag},
$$
then, $T$ is a doubly commuting $n$-tuple of regular Hilbert space oper-
ators.\par 
\end{theo3.4} 
\begin{proof} It is a consequence of Theorems 2.4 and 2.6.\par
\end{proof}

\newtheorem*{theo3.5}{Theorem 3.5}
\begin{theo3.5} Let $T=(T_1,\ldots ,T_n)$ be an $n$-tuple of commuting 
regu-
lar normal operators defined in a Hilbert space $H$. 
Next consider the $n$-tuple $b_{\beta}=(b_1^{m_1},\ldots ,b_n^{m_n})$,
where  $\beta =(m_1,\ldots ,m_n)$ is an $n$-tuple of nonnegative integers.
Then, $\prod_{i=1}^n T_i^{n_i}$ is regular
and 
$$
(\prod_{i=1}^n T_i^{n_i})^{\dag}=\prod_{i=0}^{n-1} (T_{n-i}^{\dag})^{n_i}
=\prod_{i=1}^n (T_i^{\dag})^{n_i}.
$$
\end{theo3.5}
\begin{proof}
\indent It is a consequence of Theorem 3.2.\par
\end{proof}
\indent The last application of the results of the previous section concerns Calkin 
algebras.\par
\indent Recall that if $H$ is a Hilbert space and $K(H)$ denotes the
closed ideal of all compact operators defined in $H$, then the Calkin algebra
of $H$, $C(H)=L(H)/K(H)$, is a $C^*$-algebra. Moreover, tha natural
quotient map $\pi\colon L(H)\to C(H)$ is a $C^*$-algebra morphism, see 
[5; 4] or [15; 4.1.16]. \par

\indent Furthermore, note that if $T\in L(H)$ is a regular operator,
then $\pi (T)$ is a regular element of $C(H)$. In addition, it is not difficult
to prove that if $T^{\dag}$ is the Moore-Penrose inverse of a regular
operator $T$, then $\pi (T)^{\dag}=\pi (T^{\dag})$.\par
 
\indent On the other hand, recall that an $n$-tuple $T=(T_1,\ldots ,T_n)$ is said
to be \it almost commuting (resp. almost doubly commuting)\rm, 
if $\pi (T)=(\pi (T_1),\ldots ,\pi (T_n))$ is a 
commuting tuple (resp. a doubly commuting tuple) in the $C^*$-algebra $C(H)$, 
that is if $T_iT_j-T_jT_i$ (resp. $T_iT_j-T_jT_i$ 
and $T_iT_j^*-T_j^*T_i$) belong to $K(H)$, for $i$, $j=1,\ldots ,n$, $i\ne j$,
see [6], [7] and [8]. In the following theorems the relationship between 
$n$-tuples of almost doubly commuting regular operators and 
the Moore-Penrose inverse
will be studied.\par

\newtheorem*{theo3.6}{Theorem 3.6}
\begin{theo3.6} Let $T=(T_1,\ldots ,T_n)$ be an $n$-tuple of 
regular operators defined in the Hilbert space $H$. Suppose that 
$\prod_{i=1}^nT_i$ is regular. Then, if $T$ is an almost doubly
commuting tuple of operators, 
$$
(\prod_{i=1}^nT_i)^{\dag} -\prod_{i=1}^nT_i^{\dag}\in K(H).
$$ 
\indent Conversely, if for all $i$, $j=1,\ldots ,n$, 
$i\neq j$, $T_iT_j$, $T_jT_i$, $T_i^*T_j$ and $T_jT_i^*$ are regular operators
such that
$$
(T_iT_j)^{\dag} -(T_jT_i)^{\dag}\in K(H),\hskip1cm (T_i^*T_j)^{\dag}-(T_jT_i^*)^{\dag}\in K(H),
$$
then $T$ is an almost doubly commuting $n$-tuple of regular Hilbert space operators.\par 
\end{theo3.6} 
\begin{proof} It is a consequence of Theorem 3.4 and the above recalled facts regarding
Calkin algebras.\par
\end{proof}
\newtheorem*{theo3.7}{Theorem 3.7}
\begin{theo3.7} Let $T=(T_1,\ldots ,T_n)$ be an $n$-tuple of almost commuting 
regular normal operators defined in a Hilbert space $H$. 
Next consider the $n$-tuple $T_{\beta}=(T_1^{m_1},\ldots ,T_n^{m_n})$,
where  $\beta =(m_1,\ldots ,m_n)$ is an $n$-tuple of nonnegative integers,
and suppose that $\prod_{i=1}^n T_i^{n_i}$ is regular .
Then, $T_{\beta}$ is an $n$-tuple of almost doubly commuting regular operators
defined in $H$, and 
$$
(\prod_{i=1}^n T_i^{n_i})^{\dag}-\prod_{i=1}^n (T_i^{\dag})^{n_i}\in K(H).
$$
\end{theo3.7}
\begin{proof}
\indent It is a consequence of Theorems 3.2, 3.5 and the recalled facts regarding
Calkin algebras.\par
\end{proof}
\newtheorem*{Rem3.8}{Remark 3.8}
\begin{Rem3.8}\rm It iw worth noticing that Theorems 3.2 and 3.4-3.7 can be extended
to other tuples of $C^*$-algebra elements and of Hilbert space operators respectively in the same way that
it was done in Remarks 2.3 and 2.5, that is considering permutations and
new tuples defined by Moore-Penrose inverses and  
adjoints.\par
\end{Rem3.8} 

\vskip.5truecm

\noindent Enrico Boasso\par
\noindent E-mail address: enrico\_odisseo@yahoo.it
\end{document}